\documentclass[11pt]{article}
\usepackage{amsmath}
\usepackage{graphicx}
\usepackage{amsfonts}
\usepackage{amssymb}
\newcommand{\ds}{\displaystyle}

\begin{document}

\title{{\itshape Color Visualization of Blaschke Product Mappings}}

\author{Cristina Ballantine and Dorin Ghisa}

\maketitle
 
\hfill To Professor Cabiria Andreian Cazacu\vspace{.2in}
 
\begin{abstract}
A visualization of Blaschke product mappings can be
obtained by treating them as canonical projections of covering Riemann
surfaces and finding fundamental domains and covering transformations
corresponding to these surfaces. A working tool is the technique of
simultaneous continuation we introduced  in previous papers. Here, we are
refining  this technique for some particular types of Blaschke
products,  for which coloring pre-images of annuli centered at the origin
allow us to describe the mappings with a high degree of fidelity. Additional graphics and animations are provided on the web site of the project \cite{cri}.

\end{abstract}

\section{Blaschke Products}

The building blocks of Blaschke products are M\"{o}bius transformations of
the form
\begin{equation} b_{k}(z)=e^{i\theta _{k}}\frac{a_{k}-z}{1-\overline{a}_{k}z}, \label{eq1} \end{equation}
 where  $a_{k}\in D:=\{z\in \mathbb{C}|$ $|z|<1\},$ and $\theta_{k}\in \mathbb R.$ A finite
(infinite) Blaschke product has the form
\begin{equation}w=B(z)=\prod
\limits_{k=1}^{n}b_{k}(z), \label{eq2} \end{equation}
where  $n\in \mathbb{N}$, (respectively $n=\infty).$ In the infinite case it is
customary to take \ $\displaystyle e^{i\theta_{k}}=\frac{\overline{a}_{k}}{|a_{k}|}$. It
is known that if $\sum \limits_{k=1}^{\infty}(1-|a_{k}|)<\infty,$ then the
infinite product converges uniformly on compact subsets of $W=\widehat{\mathbb C}\, 
\backslash\, (A\cup E),$ where  $A=\{1/\overline{a}_{k}\}$ and  $E$ is
the set of cluster points of  $\{a_{k}\}.$ For the sake of uniformity, we
will always take $\displaystyle \frac{\overline{a}_{k}%
}{|a_{k}|}$ for the value of $\displaystyle e^{i\theta_{k}}$ regardless whether  $B$ is finite or not.\vspace{.05in}
 
In \cite{[2]} and \cite{[3]} we gave a complete description of the domains $\Omega_{k}$ of
injectivity (fundamental domains) of the mappings when the Blaschke product
is of the form:\begin{equation} w=B_{a}(z)=\left[\frac{\overline{a}}{|a|}%
\frac{a-z}{1-\overline{a}z}\right]^{n} \label{eq3} \end{equation}
and proved the following theorem. \bigskip

\textbf{Theorem 1.1.}  \textit{The domains }$\Omega _{k}$\textit{\ bounded
by consecutive arcs of circle of the form:}
\begin{equation} z_{k}(t)=\frac{\omega _{k}t-r}{\omega _{k}tr-1}
e^{i\theta }, \ \ \ \ t\geq 0, \ \ \ k=0,1,2,...,n-1, \label{eq4} \end{equation}
\textit{where }$\omega _{k}$\textit{\ are the $n^{th}$ roots of  unity, \ }$t\geq 0$\textit{\ and }$a=re^{i\theta },$ \textit{are
fundamental domains of the covering Riemann surface }$(\widehat{\mathbb C},B_{a}).$

\vspace{.1in}

Actually, it can be easily seen that $z_{0}(t)$ is always the part of the
line (generalized circle) determined by $a$ and $1/\overline{a}$ from which
the segment between $a$ and $1/\overline{a}$ has been removed. Also, if $n$
is even, then $ \displaystyle z_{\frac{n}{2}}(t)$ is the segment between $a$ and $1/%
\overline{a}$ and the arcs $z_{\frac{n}{2}+k}(t)$ are complementary to the
arcs $z_{k}(t),$ $k=0,1,...,\frac{n}{2}-1.$ All\ the arcs $(4)$ have the end
points in  $a$ and $1/\overline{a}.$\vspace{.1in}

If we let $a\rightarrow 0,$ then $1/\overline{a} \rightarrow \infty $ \ and all the arcs $(4)$
become rays starting at the origin. The covering Riemann surface $(\widehat{\mathbb C%
},B_{0})$ is the well known Riemann surface of the multivalued function $\displaystyle
w\rightarrow w^{1/n}.$ We notice that this surface has two branch points: $0$ and $%
\infty ,$ which correspond to $a$ and $1/\overline{a}$ , the branch points
of $(\widehat{\mathbb C},B_{a}).$\ This last surface is the Riemann surface of the
multivalued function$$ w\rightarrow e^{i\theta }\frac{r-w^{1/n}}{1-rw^{1/n}}$$
obtained when trying to find an inverse of the function $(3).$\vspace{.1in}

The arcs $(4)$ have been obtained by simultaneous continuation (see \cite{[2]}) of
paths in $(\widehat{\mathbb C},B_{a})$ over the real non negative half-axis in the $%
w $-plane. We have also found\ (see \cite{[3]}) that the cover transformations of $(%
\widehat{\mathbb C},B_{a})$ were M\"{o}bius transformations of the form:\begin{equation} T_{k}(z)=\frac{a(1-\omega _{k})-(|a|^{2}-\omega _{k})z}{%
1-|a|^{2}\omega _{k}-\overline{a}(1-\omega _{k})z},\qquad k=0,1,2,...,n-1.  \label{eq5}\end{equation}
These transformations form a cyclic group of order $n$ with respect to composition. In fact, the composition rule is $T_{k}\circ T_{j}=T_{k+j\pmod n },$ $T_{0}$ being the identity and $T_{k}^{-1}=T_{n-k}.$ Every
domain $\Omega _{j}$ is mapped conformally by $T_{k}$ on the domain $\Omega
_{j+k \pmod n}$. \vspace{.1in}

The mapping $(3)$ can be visualized by drawing the pre-image of a family of
circles $w=\rho e^{i\varphi },$ $\varphi \in \lbrack 0,2\pi ),$ where $\rho $
takes different constant values. This comes to solving the equation 
$B(z)=\rho e^{i\varphi }$. The solutions are:
\begin{equation} z_{k}(\rho ,\varphi )=e^{i\theta }[r-\rho ^{1/n}e^{i(\varphi +2k\pi)/n}]/[1-r\rho ^{1/n}e^{i(\varphi +2k\pi /n)}], \ \ k=0,1,...,n-1. \label{eq6} \end{equation}

We notice that, for fixed $\rho $, every $z_{k}(\rho ,\varphi )$ represents
an arc of circle and that for every $k,$ $z_{k}(\rho ,2\pi )=z_{k+1\pmod n}(\rho ,0).$ In particular, $z_{n-1}(\rho ,2\pi )=z_{0}(\rho ,0)$ and
therefore the union of all these arcs is a full circle, which is mapped by $%
B $ on the circle $w=\rho e^{i\varphi }.$ The mapping is $n$ to $1$, every
one of the previous arcs being mapped bijectively on the circle $w=\rho
e^{i\varphi }.$ These pre-images are orthogonal to the arcs $(4)$ and have
their centers on the line determined by $a$ and $1/\overline{a}.$ An
elementary computation shows that when $\rho =1/r^{n}$ the pre-image of the
corresponding circle is the line (generalized circle) passing through $\frac{%
1}{2}(a+1/\overline{a})$ perpendicular to the line determined by $a$ and $1/%
\overline{a}.$ All the pre-image circles corresponding to $\rho <1/r^{n}$ are
on the side of the origin and those corresponding to $\rho >1/r^{n}$ are on
the other side of that perpendicular line. As $\ \rho \rightarrow 0$ these circles
accumulate to $a$ and as $\rho \rightarrow \infty $ they accumulate to $1/\overline{a}.$ We obtain an almost perfect visualization of the mapping $(3)$ by
coloring a set of annuli centered at the
origin of the $w$-plane in different colors and with saturation increasing counterclockwise and brightness increasing outward (the saturation is determined by the argument of the point and the brightness is determined by the modulus) and imposing the same color, saturation and  brightness  to the
pre-image of every point in these annuli. The Mathematica program allows the
superposition  of an orthogonal mesh on the picture in the $w$%
-plane. The mesh consists of   circles centered at the
origin and  rays starting at the origin. The pre-image of this mesh is  an orthogonal mesh in the $z$-plane.
The colors help  identify the image under $B$ of every eye of this last
mesh. Then, zooming in on the eye if necessary, we can find  the corresponding points under $B$ with high
accuracy.\vspace{.1in}

A variant of this ideea related to the case $n=1$ appears in the paper%
\textit{\ M\"{o}bius Transformations Revealed }by D.N. Arnold and J. Rogness,
Notices of AMS, Vol 55, Number 10 from November 10, 2008. Our independent
work has been presented at the International Conference on Complex Analysis
and Related Topics, Alba Iulia, Romania, August 14-19, 2008.\vspace{.1in}

 Figure 1(f) shows the fundamental domains of the mapping realized by $B_{a}$ for $a=1/2+1/3i$ (\emph{i.e.}, $ a = r e^{i\theta} \approx 0.6 e^{.59 i}$) and $n=6.$  We notice  two unbounded fundamental
domains bordered by the line $z(t) =te^{.59i}$, the two larger arcs of circle
passing through $a$ and $1/\overline{a}$, and four bounded domains bordered by the four arcs
and the segment between $a$ and $1/\overline{a}$. Each one of these six domains are mapped
conformally by B on the w-plane from which the real half axis has been removed. Figures 1(a) and (b) show pre-images of annuli under $B_{a}$. Figure 1(c) zooms in on Figure 1(a) to show the behavior  of $B_a$ close to $a$ and $1/\overline{a}$. The annuli are rendered in Figures 1 (d) and (e). Note that in all figures we show only a selection of the annuli whose pre-images are displayed. A complete collection of annuli can be viewed on the website of the project \cite{cri}.  In this example we have $1/r^{6}\approx
21.24$. Thus,  the pre-images of circles centered at the origin and of
radius less than $1/r^6$ are circles centered on the line through $a$ and $1/%
\overline{a}$ situated on the side of the origin, while the pre-images of
circles of radius greater than $1/r^6$ are situated on the other side of that
line. Taking values of $\rho $ between $0$ and $30$, we notice that indeed,
in the neighborhood of $\rho =21.24$ there is a change of convexity for the
pre-images of those circles. In order to illustrate the behavior of $B$ in
the neighborhood of $1/\overline{a}$ we have drawn  pre-images of 
annuli of very large radii. At \cite{cri} we provide an animation describing the
fundamental domains of $B_{a}$ depending on $a$  and allowing a
visualization of this dependence when  $a$ is moving inside the unit disc.\vspace{.1in}

 The next section treats a slightly different case.

\section{The Case of Two Zeros of the Same Multiplicity}

We study the case of
\begin{equation}B(z)=b_{1}^{n}(z)b_{2}^{n}(z), \mbox{\ where\ }b_{j}(z)=\frac{\overline{a}_{j}}{|a_{j}|}\frac{a_{j}-z}{1-\overline{a}_{j}z}, \ \ \ \  j=1,2\label{eq7} \end{equation}
All solitions of the equation  $B'(z)=0$ are  of order $n-1$. The solutions inside the unit circle are  $a_{1}$, $a_{2}$ and 
$b=re^{i\theta }=(1/a)[1-r_{1}^{2}r_{2}^{2}-\sqrt{%
(1-r_{1}^{2}r_{2}^{2})^{2}-|a|^{2}}]$, where $a=a_{1}a_{2}(\overline{a}_{1}+\overline{a}_{2})-(a_{1}+a_{2})$. The solution  of $B'(z)=0$  outside the unit circle is $1/\overline{b}$.

Let $\beta $ be the argument of $B(b)$ and let us perform simultaneous
continuation over the ray passing through $\displaystyle e^{i\beta }$ from all the
solutions of the equation $\displaystyle B(z)=e^{i\beta },$ \emph{i.e.},  let us solve the equation $%
B(z)=e^{i\beta }t^{n}$ for every $t>0.$ \ If we let $\lambda =e^{i\beta
/n},$ it can be easily found (see also \cite{[2]}) that:\bigskip

\textbf{Theorem 2.1.}\textit{\ The domains bounded by the curves:}\vspace{.05in}

\noindent $\displaystyle \mathit{z}_{1,2}^{(k)}\mathit{(t)=}$
\begin{equation} \mathit{\bigskip }[(r_{1}-\omega _{k}r_{2}\lambda
t)e^{i\theta _{1}}+(r_{2}-\omega _{k}r_{1}\lambda t)e^{i\theta _{2}}\pm
e^{i(\theta _{1}+\theta _{2})}\sqrt{\Delta _{k}(t)}]/2(1-\omega _{k}r_{1}r_{2}\lambda t),\label{eq8}\end{equation}
\textit{where}$$\mathit{\Delta }_{k}\mathit{(t)=[(\omega }_{k}\mathit{r}_{1}\mathit{\lambda
t-r}_{2}\mathit{)e}^{-i\theta _{1}}\mathit{-(\omega }_{k}\mathit{r}_{2}%
\mathit{\lambda t-r}_{1}\mathit{)e}^{-i\theta _{2}}\mathit{]}^{2}\mathit{%
+\bigskip }4\omega _{k}(1-r_{1}^{2})(1-r_{2}^{2})\lambda te^{-i(\theta
_{1}+\theta _{2})} $$ \textit{are fundamental domains of }$(\widehat{\mathbb C},B),$\textit{\ where }$B$%
\textit{ is given by }$(7).$\vspace{.1in}

We notice that for every $k=0,1,...,n-1$ we have $z_{1}^{(k)}(0)=a_{1},$ $\ds
\lim_{t\rightarrow \infty }z_{1}^{(k)}(t)=1/\overline{a}_{1},$ $z_{2}^{(k)}(0)=a_{2},$ 
$\ds \lim_{t\rightarrow \infty }z_{2}^{(k)}(t)=1/\overline{a}_{2}.$ An easy computation
shows that for $t>0$ the numbers $z_{1,2}^{(k)}(t)$ are all distinct except
for the case when $z_{1}^{(k)}(t)=b$ and, therefore, $z_{2}^{(k)}(t)=1/%
\overline{b}$. Thus, the curves $(8)$ can only have in common the
points $a_{j},$ respectively $1/\overline{a}_{j}$ \ and $b,$ respectively $1/%
\overline{b}$ and they are mapped bijectively by $B$ on the infinite ray
passing through $e^{i\beta }$. By an obvious extension of the conformal
correspondence theorem (see \cite{[6]}, p.154), the domains bounded by them are
represented conformally by $B$ on the complex plane from which the
previously mentioned ray has been removed.\bigskip

A visualization of the mapping $B$ can be obtained as in the previous section. Namely, we plot first the curves $z_{1,2}^{(k)}(t)$ in the $z$-plane, $k=0,1,...,n-1.$ Points in concentric annuli in the $w$-plane and their pre-images in the $z$-plane are colored as described in section 1.  If we denote $u_{k}=u_{k}(\rho ,\varphi )=\rho ^{1/n}\exp (i%
\frac{\varphi +2k\pi }{n}),$ then the equation $B(z)=\rho e^{i\varphi }$ is
equivalent to the following set of equations: \begin{equation} (1-u_{k}r_{1}r_{2})e^{-i(\theta _{1}+\theta
_{2})}z^{2}+[(r_{1}u_{k}-r_{2})e^{-i\theta
_{1}}+(r_{2}u_{k}-r_{1})e^{-i\theta
_{2}}]z+r_{1}r_{2}-u_{k}=0,\label{eq9} \end{equation}
$k=0,1,...,n-1,$
whose solutions $z_{1,2}^{(k)}(\rho ,\varphi ),$ $\varphi \in \lbrack 0,2\pi
)$ represent for every $\rho >0$ the pre-image of the circle $w=\rho
e^{i\varphi },$ $\varphi \in \lbrack 0,2\pi ).$ For $\rho \in \lbrack \rho
_{1},\rho _{2}],$ the points $z_{1,2}^{(k)}(\rho ,\varphi )$ describe the
pre-image of the annulus centered at the origin of the $w$-plane having the
radii $\rho _{1}$ and $\rho _{2}.$ When $\rho _{1}=0$ and $\rho _{2}$ is a
small number, the pre-image of the corresponding disc will be two disc-like
shaped domains around $a_{1}$ and $a_{2}.$ When $\rho _{1}$ and $\rho _{2}$
are close to $1,$ the pre-image of the corresponding annulus will be an
annulus-like shaped domain close to the unit circle.\vspace{.1in}

The intersection of these pre-images with every fundamental domain are
mapped conformally by $B$ on the previously mentioned annuli. The
fundamental domains are mapped conformally one on each other by the group of
covering transformations of $(\widehat{\mathbb C},B).$ To find these transformations
we need to solve  the equation $B(\zeta )=B(z)$ in terms of $\zeta $ (see
\cite{[5]}). This equation is equivalent to $$\frac{a_{1}-\zeta }{1-\overline{a}_{1}\zeta }\frac{%
a_{2}-\zeta }{1-\overline{a}_{2}\zeta }=\frac{a_{1}-z}{1-\overline{a}_{1}z}%
\frac{a_{2}-z}{1-\overline{a}_{2}z}\omega _{k},\qquad k=0,1,...,n-1.$$
Each one of these second degree equations has two solutions $\zeta =S_{k}(z)$
and $\zeta ^{\prime }=T_{k}(z)$ and the relation between $\zeta $ and $\zeta
^{\prime }$ is given by:$$\frac{a_{1}-\zeta }{1-\overline{a}_{1}\zeta }\frac{%
a_{2}-\zeta }{1-\overline{a}_{2}\zeta }=\frac{a_{1}-\zeta ^{\prime }}{1-%
\overline{a}_{1}\zeta ^{\prime }}\frac{a_{2}-\zeta ^{\prime }}{1-\overline{a}%
_{2}\zeta ^{\prime }},$$
which is equivalent to 
$$\zeta ^{\prime }=U(\zeta )=\frac{A-\zeta }{1-\overline{A}\zeta }, \mbox{\ \ 
where\ \ } A=\frac{a_{1}a_{2}(\overline{a}_{1}+\overline{a}%
_{2})-(a_{1}+a_{2})}{|a_{1}a_{2}|^{2}-1}$$
We notice that $U$ does not depend on $k$ and it is an involution. Moreover, 
$U\circ S_{k}=T_{k}$ and $U\circ T_{k}=S_{k}$ for every $k$ and $%
U(a_{1})=a_{2},$ $U(1/\overline{a}_{1})=1/\overline{a}_{2},$ $%
U(a_{2})=a_{1}, $ $U(1/\overline{a}_{2})=1/\overline{a}_{1}.$ In other words,
we have proved:\bigskip

\textbf{Theorem 2.2}. \textit{The group of covering transformations of }$(%
\widehat{\mathbb C},B)$\textit{\ is the group generated by }$S_{k},$\textit{\ }$%
k=0,1,...,n-1,$\textit{\ and }$U.$\textit{\ The composition rule is }$U\circ
S_{k}=T_{k},$\textit{\ }$U\circ T_{k}=S_{k},$\textit{\ }$S_{k}\circ
S_{j}=S_{k+j\pmod n},$\textit{\ particularly }$S_{k}^{-1}=S_{n-k},$%
\textit{\ }$T_{k}^{-1}=(U\circ S_{k})^{-1}=S_{k}^{-1}\circ U=S_{n-k}\circ U$%
\textit{. The identity of the group is }$S_{0}$\textit{.}

\ \ 

The solutions of the equations $(9)$
are\bigskip\ functions of $u_{k}=u_{k}(\rho ,\varphi )$ of the form:
\begin{equation}z_{1,2}^{(k)}(\rho ,\varphi )=[(r_{1}-r_{2}u_{k})e^{i\alpha
_{1}}+(r_{2}-r_{1}u_{k})e^{i\alpha _{2}}\pm \sqrt{\Delta _{k}}%
]/2(1-r_{1}r_{2}u_{k}),\label{eq11}\end{equation}
where $$\Delta _{k}=\Delta _{k}(\rho ,\varphi )=[(r_{1}-r_{2}u_{k})e^{i\alpha
_{1}}-(r_{2}-r_{1}u_{k})e^{i\alpha
_{2}}]^{2}+4(1-r_{1}^{2})(1-r_{2}^{2})u_{k}e^{i(\alpha _{1}+\alpha _{2})}.$$

We notice that when $r_{2}\rightarrow1,$ then
$\Delta _{k}(\rho ,\varphi )\rightarrow[(r_{1}-u_{k})e^{i\alpha
_{1}}-(1-r_{1}u_{k})e^{i\alpha _{2}}]^{2}$ and therefore $z_{1}^{(k)}(\rho ,\varphi )\rightarrow e^{i\alpha _{1}}[r_{1}-u_{k}(\rho
,\varphi )]/[1-r_{1}u_{k}(\rho ,\varphi )],$
which are, for constant $\rho ,$ circles orthogonal to the boundaries of the
fundamental domains of $ [b_{1}(z)]^{n}.$ At the same time, $%
z_{2}^{(k)}(\varphi )\rightarrow e^{i\alpha _{2}}$ for every $k=0,1,...,n-1.$\bigskip

To visualize the kind of mapping a product $(7)$ produces, we have chosen
again $n=6$ , we took $a_{1}$ to be the value $a$ from the Figure 1 and allowed
$a_{2}$ to take different values in the unit disk. An animation provided on the web
site \cite{cri} shows how the fundamental domains of $B$ from Figure 1(f) change when
introducing $a_{2}$ in the picture. In particular we notice that, as $a_{2}$ approaches a
point $e^{i\alpha _{2}}$, some of these domains shrink  to $%
e^{i\alpha _{2}},$ while the others take the shape of those corresponding to
the Blaschke product $b_{1}^{n}.$  We notice that those changes are
important only when $a_{2}$ is introduced in an unbounded fundamental domain
of $B.$ When $a_{2}$ belongs to a bounded fundamental domain of $B,$ they
are not so important and if $a_{2}$ is close to the unit circle in such a
domain, the changes are hardly noticeable. Such a situation happens usually
when, for an infinite Blaschke product $B,$ we are switching from a partial
product $B_{n}$ of $B$ to the next partial product $B_{n+1}.$ This remark
suggests that Blaschke factors generated by zeros close to the unit circle
do not influence too much the picture, and therefore not only $(B_{n})$
approximate uniformly $B$ on compact subsets of $\widehat{\mathbb C}$ $\backslash $ $%
(A\cup E),$ but the fundamental domains of $(B_{n})$ "approximate" in turn
the fundamental domains of $B.$

\vspace{.1in}

Figure 2 illustrates the Blaschke Product $(7)$ with two zeros of order $6$.  We have taken $\ds a_1 = \frac{1}{2} +\frac{1}{3} i$, $a_2 = \frac{4}{5}e^{2.5 i}$ and $n=6$. Figure 2(f) shows the fundamental domains of $B$. Eleven of the fundamental  are bounded and only one is unbounded. The bounded domains  are
bordered by arcs ending in $a_1$ and $1/\overline{a}_1$, $a_2$ and $1/\overline{a}_2$ as well as $b$ and
$1/\overline{b}$, while the unbounded domain is bordered by the arcs ending in $a_1$ and
$1/\overline{a}_1$ and $b$ and $1/\overline{b}$. The Blaschke product $B$ mapps conformally each one of
these domains on the w-plane from which the ray passing through $B(b)$ has been
removed. Figures 2(a) and (b) show the pre-images of concentric annuli. Some of the annuli are shown in Figure 2(d) and (e). A complete collection of the annuli is displayed at the web site \cite{cri}. Figure 2(c) zooms in on the pre-image of $B$ close to $a_2$ and $1/\overline{a}_2$. Figures 2(g) and (h) show the fundamental domains of $B$ above with $a_2$ replaced by  $.99e^{2.5 i}$ and $ .999 e^{2.5 i}$ respectively. 

\section{The Case of Zeros of the Same Module and
of Arguments $\alpha +2k\pi /n$}

\bigskip

Suppose that  $a_{k}=re^{i\alpha}\omega_{k}$, $k=0,1,2,...,n-1$,
where \ $\omega_{k}$ are the roots of order $n$ of  unity. Then  $\ds \frac {\overline{a}_{k}}{|a_{k}|}=e^{-i\alpha}\omega_{n-k}$. Therefore $\ds \prod\limits_{k=0}^{n-1}\frac{\overline{a}_{k}}{|a_{k}|}=(-1)^{n-1}e^{-n\alpha i}$
and
\begin{equation}w=B(z)=e^{-in\alpha }\frac{z^{n}-r^{n}e^{in\alpha }}{e^{-in\alpha }r^{n}z^{n}-1}\label{eq13}\end{equation}\vspace{.1in}

\textbf{Theorem 3.1}. \textit{The\ domains bounded by consecutive rays }$z_{k}(t)=e^{i[\alpha +(2k+1)\pi /n]}t,$\textit{\ }$t\geq 0$\textit{\ are
fundamental domains of }$\ (\widehat{\mathbb C},B),$\textit{\ where }$B$\textit{\ is
given by }$(11).$\textit{\ The image by }$B$\textit{\ of every one of these
domains is the }$w$\textit{-plane from which the interval }$(r^{n},1/r^{n})$%
\textit{\ of the real axis has been removed. The covering transformations of 
}$(\widehat{\mathbb C},B)$\textit{\ are rotations }$T_{j}(z)=\omega _{j}z,$\textit{\
where }$\omega _{j}$\textit{\ are the roots of order }$n$\textit{\ of 
unity}$.$

\vspace{.1in}

\textbf{Proof:} The equation $B(z)=t,$ $t\geq 0$ (simultaneous continuation
over the real non negative half-axis) has the solutions 
\begin{equation}z_{k}(t)=\left\{
\begin{array}{ll} [(r^{n}-t)/(1-r^{n}t)]^{1/n}e^{i\alpha }\omega _{k}& \ \ \ \ \mbox{if} \ t\in \lbrack
0,r^{n}]\cup (1/r^{n},\infty ) \\ 
\\ 
\lbrack (t-r^{n})\text{ }/\text{ }(1-r^{n}t)]^{1/n}e^{i(\alpha +\pi
/n)}\omega _{k}&\ \ \ \ \mbox{if} \ t\in \lbrack r^{n},1/r^{n}),%
\end{array}\right.
\label{eq14}\end{equation}
where $k=0,1,...,n-1$\bigskip

We notice that $z_{k}(r^{n})=0$ and $\ds\lim_{t\rightarrow1/r^{n}}z_{k}(t)=\infty $ 
for every $k,$ and as $t$ varies between $r^{n}$ and  $1/r^{n}$ the
argument of $z_{k}(t)$ remains the same, namely $\ds \alpha +(2k+1)\frac{\pi }{n}.$ Therefore, $z_{k}(t)$ describes the ray issued from the origin and
forming  the angle $\alpha +(2k+1)\frac{\pi }{n}$ with the positive real half-axis. The intersection of this ray with the unit circle is the point $\ds \zeta
_{k}=z_{k}(1)=e^{[\alpha +(2k+1)\pi /n]i}$ and  it obviously does not depend
on $r.$ This means that the respective ray will remain the same as  $r$ 
varies between  $0$  and  $1$. 

\vspace{.1in}

The above mentioned rays are the borders of fundamental domains $\Omega
_{k},$ which are mapped conformally by $B$ on the domain obtained when
removing from the $w$-plane the interval $[r^{n},1/r^{n}].$ The expressions
of $z_{k}(t),$ for $t\in \lbrack 0,r^{n}]\cup (1/r^{n},\infty )$ show simply
that the image by $B$ of every segment between $0$ and $a_{k}$ is the
segment $[0,r^{n}]$, with $B(0)=$ $r^{n},$ $B(a_{k})=0,$ and that of every
interval $(1/a_{k},\infty )$, on the ray $z(t)=e^{i\alpha }\omega _{k}t$, is
the interval $(1/r^{n},\infty ),$ with $B(\infty )=1/r^{n}$ and $%
B(1/a_{k})=\infty $ for every $k=0,1,...,n-1.$ A small circle around the
origin of the $w$-plane has as pre-image by $B$ the union of $n$ curves,
each one around an $a_{k}$ and situated in $\Omega _{k}.$ Every one of those
curves is mapped by $B$ bijectively on the respective circle. A very big
circle around the origin\ of the $w$-plane has as pre-image by $B$ the union
of $n$ small curves each one around a point $1/\overline{a}_{k}.$ Every one
of these curves is mapped by $B$ bijectively on the respective circle. A
circle $w(\tau )=\rho e^{i\tau },$ $\tau \in \lbrack 0,2\pi ),$ with $\rho $
close to $1$ has as pre-image by $B$ a unique curve close to the unit
circle. Its intersection  with every $\Omega _{k},$ is mapped by $B$ bijectively on the respective circle. That curve is orthogonal to every
ray $\ds z(t)=e^{i[\alpha +(2k+1)\pi /n]}t$,  $t>0$. \vspace{.1in}

The domains $\Omega _{k}$ are mapped one on each other by the covering
transformations $T_{j}(z)=\omega _{j}z.$ It can be easily checked that,
for every $j$, we have $B\circ T_{j}=B$, \emph{i.e.}, $T_{j}$ are indeed covering
transformations. Obviously, the set of transformations $\{T_{j}\}$ is a
cyclic group of order $n$ for which the composition law is $T_{j}\circ
T_{k}=T_{j+k \pmod n}.$ Every  $\Omega _{k}$ is mapped conformally by
 $T_{j}$ on $\Omega _{k+j \pmod n}.$

\vspace{.1in}

In  Figure 3 we consider the Blaschke product  of oder $6$  defined in $(11)$ with zeroes $a_k=re^{i \alpha}\omega_k$, where $r=2/3$, $\alpha=\pi/5$ and $k=0,1,\ldots, 5$. Figure 3(g) shows the fundamental domains of $B$ and  represents six infinite sectors bounded by the rays \ $t\rightarrow e^{i[\alpha +(2k+1)\pi /6]},k=0,1,\dots,5$. Every such sector is mapped by $B$ conformally on  the $w$-plane from which the
interval $(r^{6},1/r^{6})$ has been removed. Figures 3(a) and (b) show the pre-images under $B$ of concentric annuli. Several of these annuli have been displayed in figures 3(e) and (f). The complete collection of annuli is available at \cite{cri}. Every small annulus around the
origin in the $w$-plane has as pre-image six annular domains, each one
around an $a_{k}.$ The spectrum of colors in every one of these domains
coincide with that of the annulus. Every annulus of radii close to $1$ has
as pre-image a unique annular domain close to the unit circle. The spectrum
of colors of the intersection of this annular domain with every $\Omega _{k}$
coincide with that of the annulus. Figure 3(c) zooms in on the pre-image of $B$ close to one of the zeros.  Figure 3(d) zooms in on the pre-image of $B$ close to the origin.

\section{Two Sets of Zeros of Arguments $\alpha _{1}+2k\pi /n$
and $\alpha _{2}+2k\pi /n$}

Following the pattern of the section 2, for $0\leq r_{1}\leq r_{2}<1$,
we consider  Blaschke products having the zeros  $a_{k}=r_{1}e^{i\alpha
_{1}}\omega _{k}$, and  $b_{k}=r_{2}e^{i\alpha _{2}}\omega _{k}$,
namely

\begin{equation}B(z)=\frac{(z^{n}e^{-in\alpha _{1}}-r_{1}^{n})(z^{n}e^{-in\alpha
_{2}}-r_{2}^{n})}{(r_{1}^{n}e^{-in\alpha
_{1}}z^{n}-1)(r_{2}^{n}e^{-in\alpha _{2}}z^{n}-1)}\label{eq15}\end{equation}

In the following we will study in detail the case $\alpha _{1}=\alpha
_{2}=\alpha $ and at the end of the section we will make some remarks about the general
case.

\vspace{.1in}

\textbf{Theorem 4.1.}\textit{\ The covering Riemann surface }$(\widehat{\mathbb C}%
,B),$\textit{\ where }$B$\textit{\ is given by }$(13)$, \textit{\ has as
fundamental domains the sectors bounded by consecutive infinite rays }$\ds
z(\tau )=e^{i(\alpha +k\pi /n)}\tau ,$\textit{\ }$\tau \geq 0,$\textit{\ }$%
k=0,1,...,2n-1.$\textit{\ The group of covering transformations is generated
by the uniform branches of the multivalued function}$$z\rightarrow
\{[(1+r_{1}^{n}r_{2}^{n})z^{n}-(r_{1}^{n}+r_{2}^{n})e^{in\alpha
}]/[(r_{1}^{n}+r_{2}^{n})e^{-in\alpha }z^{n}-(1+r_{1}^{n}r_{2}^{n})]\}^{1/n}$$
\textit{and the transformation }$z\rightarrow e^{2\pi i/n}z.$ \textit{ The
respective sectors are mapped conformally by }$B$\textit{\ on\ the }$w$%
\textit{-plane from which a part of real axis has been removed}$.$\vspace{.1in}

\textbf{Proof:} The equation $B(z)=t,$ where $B$ is given by $(13)$  has
$2n$ solutions, which can be written in the form:
\begin{equation} z_{1,2}^{(k)}(t)=\left\{ \begin{array}{lr}e^{i\alpha }\omega _{k}K_{1,2}^{1/n}(t) &\ \ \ \mbox{ if}\ \  K_{1,2}(t)>0 \\ \ \\ e^{i(\alpha +\pi /n)}\omega _{k}[-K_{1,2}(t)]^{1/n} &\ \ \ \mbox{if}\ \  K_{1,2}(t)<0 \end{array}\right. \label{eq16} \end{equation}
where
\begin{equation} K_{1,2}(t):=\frac{(r_{1}^{n}+r_{2}^{n})(1-t)\pm \sqrt{\Delta (t)}}{2(1-r_{1}^{n}r_{2}^{n}t)}\label{eq17} \end{equation}
and
\begin{equation}\Delta
(t)=(r_{1}^{n}+r_{2}^{n})^{2}(1-t)^{2}-4(r_{1}^{n}r_{2}^{n}-t)(1-r_{1}^{n}r_{2}^{n}t). 
\label{eq18} \end{equation}\vspace{.1in}

We notice that $\sqrt{\Delta (1)}=2(1-r_{1}^{n}r_{2}^{n})$ and  hence $
K_{1}(1)=1$ and $K_{2}(1)=-1$, which implies that  $z_{1}^{(k)}(1)=e^{i%
\alpha }\omega _{k}$ and  $z_{2}^{(k)}(1)=e^{i(\alpha +\pi /n)}\omega
_{k}. $ Indeed, it is  easy to check that $B(e^{i\alpha }\omega
_{k})=B(e^{i(\alpha +\pi /n)}\omega _{k})=1.$ On the other hand, we have
$$z_{2}^{(k)}(r_{1}^{n}r_{2}^{n})=0 \ \ \ \mbox{and} \ \ \ 
z_{1}^{(k)}(r_{1}^{n}r_{2}^{n})=\frac{e^{i\alpha }\omega
_{k}(r_{1}^{n}+r_{2}^{n})^{1/n}}{(1+r_{1}^{n}r_{2}^{n})^{1/n}}.$$

This shows that, as  $t$ varies from  $r_{1}^{n}r_{2}^{n}$ to  $1$
every point $z_{1}^{(k)}(t)$  moves on the ray $e^{i\alpha }\omega
_{k}\tau$, $\tau \geq 0$, from  $e^{i\alpha }\omega
_{k}(r_{1}^{n}+r_{2}^{n})^{1/n}/(1+r_{1}^{n}r_{2}^{n})^{1/n}$ to the
unit circle, while every point $z_{2}^{(k)}(t)$ moves from the origin to
the unit circle on the ray $e^{i(\alpha +\pi /n)}\omega _{k}\tau$,  $\tau
\geq 0$. This time the simultaneous continuation over $w(t)=t,$ $0\leq t\leq
1$, does not bring us to a satisfactory situation, since not all  corresponding
paths  meet  in  branch points of $(\widehat{\mathbb C},B).$ These branch
points can be found by solving the equation $B^{\prime }(z)=0.$ The localization problem of the zeros of B'(z)
will reoccur in the following sections (see \cite{[7]} for more details). The solutions
of the equation $B'(z)=0$ are
$c_{1,2}^{(k)}=e^{i\alpha }\omega _{k}\rho
_{1,2}, \ k=0,1,...,n-1,$
where
\begin{equation}\rho _{1,2}=\frac{(r_{1}^{n}r_{2}^{n}+1\pm \sqrt{%
(r_{1}^{n}r_{2}^{n}+1)^{2}-(r_{1}^{n}+r_{2}^{n})^{2}})^{1/n}}{
(r_{1}^{n}+r_{2}^{n})^{1/n}}.\label{eq19}\end{equation}

For every $k,$ the points $c_{1}^{(k)}$ and $c_{2}^{(k)}$ are symmetric with
respect to the unit circle, since $\rho _{1}\rho _{2}=1$. We notice also
that $r_{1}<\rho _{2}<r_{2}$, hence $1/r_{2}<\rho _{1}<1/r_{1}$ and
that $-1<B(c_{2}^{(k)})<0,$ hence $B(c_{1}^{(k)})<-1.$ We conclude from here
that as $\tau $ varies between $r_{1}$and $r_{2},$ the point $B(e^{i\alpha
}\omega _{k}\tau )$ varies on the real
negative half-axis between $B(c_{2}^{(k)})$ and $0$. Similarly, as $\tau $ varies between $1/r_{2}$ and $1/r_{1},$ the point $B(e^{i\alpha }\omega _{k}\tau )$ varies on the real negative half-axis between $B(c_{1}^{(k)})$ and $-\infty $. Therefore the
continuation should also be performed  over the intervals from $%
B(c_{2}^{(k)})$ to $0,$ respectively from $B(c_{1}^{(k)})$ to $-\infty
. $ In other words, the image of any ray passing through $e^{i\alpha }\omega
_{k}$ is the set obtained by removing  the interval
from $B(c_{1}^{(k)})$ to $B(c_{2}^{(k)})$ from the real axis.\vspace{.1in} 

 Let us denote by $\Omega _{k}^{\prime }$ the domains bounded by the
rays
$$z=\tau e^{i(\alpha -\pi /n)}\omega _{k},\tau \geq 0\ \ \ \mbox{and} \ \ \ 
z=\tau e^{i\alpha }\omega _{k}, \tau \geq 0, \ \ \ k=0,1,2,...,n-1$$
and by $\Omega _{k}^{\prime \prime }$ the domains bounded by the rays
$$z=\tau e^{i\alpha }\omega _{k}, \tau \geq 0 \ \ \ \mbox{and} \ \ \ z=\tau
e^{i(\alpha +\pi /n)}\omega _{k}, \tau \geq 0,\ \ \ k=0,1,2,...,n-1.$$
By the conformal correspondence theorem (see \cite{[6]}, page 154), the domains $%
\Omega _{k}^{\prime }$ and $\Omega _{k}^{\prime \prime }$ are mapped
conformally by $B$ on the complex plane from which the part of the real axis
between $B(c_{2}^{(k)})$ and $+\infty $ and between $-\infty $ and $%
B(c_{1}^{(k)})$ has been removed.\vspace{.1in}

The website \cite{cri} provides several plots  showing how these fundamental
domains change when varying $r_{1}$ and/or $r_{2}.$ We notice that the
intersections $e^{i\alpha }\omega _{k}$ and $e^{i(\alpha +\pi /n)}\omega
_{k} $ of the above mentioned rays with the unit circle remain the same,
\emph{i.e.},  they do not depend on the particular values of $r_{1}$ and $r_{2},$ but
only on $\alpha $ and $n.$ On the other hand, if we let just one point $a_{k}$
move on the ray $e^{i\alpha }\omega _{k}\tau ,$ then all these intersections
will suffer perturbations.\vspace{.1in}

Notice also that $B(\omega _{k}z)=B(z)$ and therefore $T_{k}(z)=\omega
_{k}z$ are covering transformations of the covering Riemann surface $(%
\widehat{\mathbb C},B).$ These are not the only such transformations. In order to find all the
covering transformations of $(\widehat{\mathbb C},B)$ we need to find, as in the
previous section, $\zeta $ as a function of $z$ for which $B(\zeta )=B(z).$
It can be easily checked that this happens if and only if either $\zeta
^{n}=z^{n}$ (hence $\zeta =\omega _{k}z),$ or
 \begin{equation} \zeta ^{n}=\frac{(1+r_{1}^{n}r_{2}^{n})z^{n}-(r_{1}^{n}+r_{2}^{n})e^{in\alpha }}{(r_{1}^{n}+r_{2}^{n})e^{-in\alpha }z^{n}-(1+r_{1}^{n}r_{2}^{n}) } \label{eq20}\end{equation}

For $z=\rho e^{i(\alpha \pm \pi /n)}\omega _{k},$ $\rho \geq 0,$ we have $z^{n}=-\rho ^{n}e^{in\alpha }$ and the expression $(18)$ becomes $$\zeta ^{n}=\frac{([(1+r_{1}^{n}r_{2}^{n})\rho
^{n}+(r_{1}^{n}+r_{2}^{n})]e^{in\alpha }}{(r_{1}^{n}+r_{2}^{n})\rho
^{n}+(1+r_{1}^{n}r_{2}^{n})}$$
with the solutions \begin{equation}\zeta _{j}(\rho )=\left(\frac{(1+r_{1}^{n}r_{2}^{n})\rho ^{n}+(r_{1}^{n}+r_{2}^{n})}{(r_{1}^{n}+r_{2}^{n})\rho
^{n}+(1+r_{1}^{n}r_{2}^{n})}\right)^{1/n}e^{i\alpha }\omega _{j}, \label{21} \end{equation}
where $j=0,1,2,...,n-1$. \vspace{.1in}

For $z=\rho e^{i\alpha }\omega _{k},$ $\rho \geq 0,$ we have $z^{n}=\rho
^{n}e^{in\alpha }$ and the expression $(18)$ becomes$$\zeta ^{n}=\frac{[(1+r_{1}^{n}r_{2}^{n})\rho
^{n}-(r_{1}^{n}+r_{2}^{n})]e^{in\alpha
}}{(r_{1}^{n}+r_{2}^{n})-(1+r_{1}^{n}r_{2}^{n})} $$
 with the solutions
\begin{equation}\zeta _{j}^{\prime }(\rho )=\left(\frac{(1+r_{1}^{n}r_{2}^{n})\rho
^{n}-(r_{1}^{n}+r_{2}^{n})}{(r_{1}^{n}+r_{2}^{n})\rho
^{n}-(1+r_{1}^{n}r_{2}^{n})}\right)^{1/n}e^{i\alpha }\omega _{j}, \ \ \  j=0,1,2,...,n-1,\label{eq22} \end{equation}
if $ \rho \in (0,$\ $[r_{1}^{n}+r_{2}^{n}]^{1/n}$ $/$ $%
[1+r_{1}^{n}r_{2}^{n}]^{1/n})\cup ([1+r_{1}^{n}r_{2}^{n}]^{1/n}/$ $%
[r_{1}^{n}+r_{2}^{n}]^{1/n},+\infty )$, and
\begin{equation}\zeta _{j}^{\prime }(\rho )=\frac{-(1+r_{1}^{n}r_{2}^{n})\rho
^{n}+(r_{1}^{n}+r_{2}^{n})}{[(r_{1}^{n}+r_{2}^{n})\rho
^{n}-(1+r_{1}^{n}r_{2}^{n})]^{1/n}}e^{i(\alpha +\pi /n)}\omega _{j} \ \ \  j=0,1,2,...,n-1,\label{eq23} \end{equation}
if $\rho \in ([r_{1}^{n}+r_{2}^{n}]^{1/n}$ $/$ $%
[1+r_{1}^{n}r_{2}^{n}]^{1/n},[1+r_{1}^{n}r_{2}^{n}]^{1/n}/$ $%
[r_{1}^{n}+r_{2}^{n}]^{1/n})$. \vspace{.1in}

If we choose the uniform branches $S_{j}(z)$ of the multivalued
function
\begin{equation}\left(\frac{(1+r_{1}^{n}r_{2}^{n})z^{n}-(r_{1}^{n}+r_{2}^{n})e^{in\alpha }}{(r_{1}^{n}+r_{2}^{n})e^{-in\alpha}z^{n}-(1+r_{1}^{n}r_{2}^{n})}\right)^{1/n}\label{eq24} \end{equation}
such that $S_{j}(\rho e^{i(\alpha -\pi /n)}\omega _{k})=\zeta _{j+k}(\rho ),$
then $S_{j}$ maps conformally the interior of every $\Omega _{k}^{\prime }$
on the interior of $\Omega _{k+j}^{\prime \prime }$ and\ the interior of
every $\Omega _{k}^{\prime \prime }$ on the interior of $\Omega
_{k+j+1.}^{\prime }$ Indeed, we notice that 
\begin{equation}S_{j}(0)=\zeta _{j+k}(0)=\frac{[r_{1}^{n}+r_{2}^{n}]^{1/n}e^{i\alpha
}\omega _{j+k}}{[1+r_{1}^{n}r_{2}^{n}]^{1/n}}=\zeta _{j+k}^{\prime
}(0)\label{eq25}\end{equation}
and
\begin{equation}S_{j}(\infty )=\zeta _{j+k}(\infty
)=\frac{[1+r_{1}^{n}r_{2}^{n}]^{1/n}e^{i\alpha }\omega _{j+k}}{
r_{1}^{n}+r_{2}^{n}}=\zeta _{j+k}^{\prime }(\infty )\label{eq26}\end{equation}

As $z$ moves on the ray $z(\rho )=\rho e^{i(\alpha -\pi /n)}\omega _{k},$ $%
\rho \geq 0$ from $0$ to $\infty ,$ $S_{j}(z(\rho ))$ moves on the ray $%
z(\tau )=\tau e^{i\alpha }\omega _{j+k}$ between the two values $\zeta
_{j+k}(0)$ and $\zeta _{j+k}(\infty ).$ Also, as $z$ varies on the ray $%
z(\rho )=\rho e^{i\alpha }\omega _{k}$ from $0$ to $\zeta _{k}(0)$ and from $%
\zeta _{k}(\infty )$ to $\infty ,$ $S_{j}(z(\rho ))$ varies from $\zeta
_{j+k}(0)$ to $0,$ respectively from $\infty $ to $\zeta _{j+k}(\infty ).$
In other words, the ray $z(\tau )=\tau e^{i\alpha }\omega _{j+k},$ $\tau
\geq 0$ is the bijective image by $S_{j}$ of the path formed with the ray $%
z(\rho )=\rho e^{i(\alpha -\pi /n)}\omega _{k}$ and the two intervals from $%
0 $ to $\zeta _{k}(0)$ and from $\zeta _{k}(\infty )$ to $\infty $ of the
ray $z(\rho )=\rho e^{i\alpha }\omega _{k}.$ On the other hand, if we let $z$
vary on the ray $z(\rho )=\rho e^{i\alpha }\omega _{k}$ between $\zeta
_{k}(0)$ and $\zeta _{k}(\infty ),$ $S_{j}(z(\rho ))$ describes the ray $\ds
z(\tau )=\tau e^{i(\alpha +\pi /n)}\omega _{j+k}.$ The fact that the
mappings we mentioned above are conformal is again a corollary of the
conformal correspondence theorem. Moreover, it is obvious that every mapping 
$S_{j}$ is conformal throughout, except for the points 
$$
0, \ \ \left(\frac{r_{1}^{n}+r_{2}^{n}}{1+r_{1}^{n}r_{2}^{n}}\right)^{1/n}e^{i\alpha
}\omega _{k},\ \  \left(\frac{1+r_{1}^{n}r_{2}^{n}}{r_{1}^{n}+r_{2}^{n}}\right)^{1/n}e^{i\alpha
}\omega _{k},  \ \ \mbox{and} \ \ \infty.$$  It can be easily checked that $%
B(S_{j}(z))=B(z)$ for every $j$ and hence $S_{j}$ are covering transformations
of $(\widehat{\mathbb C},B).$\vspace{.1in}

In order to find the composition law of the mappings $S_{j},$ we notice that
the ray $\ds z(\rho )=\rho e^{i(\alpha -\pi /n)}\omega _{k}$ is mapped
bijectively by $S_{j^{\prime }}$ on the interval between $S_{j^{\prime }}(0)$
and $S_{j^{\prime }}(\infty )$ on the ray $\ds z(\tau )=\tau e^{i\alpha }\omega
_{k+j^{\prime }},$ which is mapped in turn bijectively by $S_{j}$ on the ray 
$z(t)=te^{i\alpha }\omega _{k+j^{\prime }+j}$. Therefore, $S_{j}\circ
S_{j^{\prime }}$ maps bijectively the ray $\ds z(\rho )=\rho e^{i(\alpha -\pi
/n)}\omega _{k}, \ \rho \geq 0$ on the ray $\ds z(t)=te^{i\alpha }\omega
_{k+j^{\prime }+j}$ with $S_{j}\circ S_{j^{\prime }}(0)=0$ and $S_{j}\circ
S_{j^{\prime }}(\infty )=\infty .$ In a similar way we find that the ray $\ds
z(\rho )=\rho e^{i\alpha }\omega _{k}, \ \rho \geq 0$ is mapped bijectively by $%
S_{j}\circ S_{j^{\prime }}$ on the ray $\ds z(t)=te^{i(\alpha +\pi /n)}\omega
_{k+j^{\prime }+j}, \ t\geq 0,$ which implies that the domain $\Omega
_{k}^{\prime }$ is mapped conformally by $S_{j}\circ S_{j^{\prime }}$ on the
domain $\Omega _{k+j^{\prime }+j}^{\prime }$ with $S_{j}\circ S_{j^{\prime
}}(0)=0,$ $S_{j}\circ S_{j^{\prime }}(\infty )=\infty $ and the mapping is
conformal on the closure of $\Omega _{k}^{\prime },$ except perhaps in $0$
and $\infty .$ In a similar way we find that $S_{j}\circ S_{j^{\prime }\text{
}}$maps conformally $\Omega _{k}^{\prime \prime }$ on $\Omega _{k+j^{\prime
}+k}^{\prime \prime }$ and this happens for every $k.$ It is a simple fact
that such a transformation must coincide with $z\rightarrow \omega _{j^{\prime }+j}z$. 
In other words $S_{j}\circ S_{j^{\prime }}=T_{j+j^{\prime }}.$ It is also
obvious that $S_{j}\circ T_{k}=T_{k}\circ S_{j}=S_{j+k}.$\vspace{.1in}

Figure 4 illustrates the Blaschke product defined by $(13)$ with 
 $n=4$, $\alpha=\pi/3$, $r_1=3/5$ and $r_2=4/5$. Thus, $a_k=3/5 e^{\pi/3}\omega_k$ and $b_k=4/5 e^{\pi/3}\omega_k$, $k=0,1,2,3$.  Figure 4(g) shows the fundamental domains of $B$. It represent the unit circle
and 8 infinite rays passing through $e^{i(\alpha +k\pi /4)},k=0,1,...,7$ .
The points $a_{j}$, $b_{j}$ and $c_{1,2}^{(j)}$, are on the ray $%
z\rightarrow e^{i(\alpha +j\pi /2)}z,$ $j=0,1,2,3$.  Figures 4(a) and (b) show the pre-images under $B$ of concentric annuli. Several of these annuli have been displayed in figures 4(c) and (d). The complete collection of annuli is available at \cite{cri}. Figures 4(c) and (d) zoom in on the pre-image plot close to the zeros in the second quadrant inside and outside the unit circle respectively. \vspace{.1in}

The study of the general case when $\alpha _{1}$ and $\alpha _{2}$ are
arbitrary is more laborious, yet still feasible. The images by $B$ of the
branch points are no longer real if $\alpha _{1}\neq \alpha _{2}$ and the
simultaneous continuation should be performed over a path passing through
those images. Moreover, the equation for the branch points is  more
complicated. We present below the final results of our computation. With the
notations \ $a=r_{1}e^{i\alpha _{1}}$ and \ $b=r_{2}e^{i\alpha _{2}}$ the
equation $B^{\prime }(z)=0$ becomes
\begin{equation}z^{n-1}[\overline{A}z^{2n}-2z^{n}+A]=0,\label{eq27}\end{equation}
where \begin{equation}A=\frac{a^{n}+b^{n}-a^{n}b^{n}(\overline{a}^{n}+\overline{b}^{n})}{
1-|a|^{2n}|b|^{2n}}.\label{eq28}\end{equation}
The origin is a branch point of order $n$ (and so is $z=\infty $). There are 
$2n$ more distinct branch points. If we denote by $b_{k}$ the non zero
solutions of the above equation  situated in the unit disc,\ then $1/%
\overline{b}_{k},\ k=0,1,...,n-1$ are also solutions and they are situated   in
the exterior of the unit disk. Let $\beta $ be the argument of $B(b_{k})$, 
which does not depend on $k,$ and it is the same as the argument of $B(1/%
\overline{b}_{k}).$ We need to perform simultaneous continuation over the ray $
w=e^{i\beta }t,$ for $\ t$ between $1$ and $|B(b_{k})|$ and then over the
segment between $B(b_{k})$ and $B(0).$ After reflection in the unit circle,
the curves obtained will give us the boundaries of the fundamental domains
of $B.$ It is obvious that $B(\omega _{k}z)=B(z)$ for every $n$-th root of
unity $\omega _{k},$ hence $T_{k}(z)=\omega _{k}z$ are cover transformations
of the covering Riemann surface $(\widehat{\mathbb C},B).$ On the other hand, solving
the equation $B(\zeta )=B(z)$ in terms of $\zeta $ we obtain  $\zeta
=z\omega _{k},$ $k=0,1,...,n-1,$ or
\begin{equation} \zeta ^{n}=\frac{z^{n}-A}{\overline{A}z^{n}-1},\label{eq24'}\end{equation}
where $A$ is given by $(26).$\ \ 

Let us denote by $\zeta =S_{j}(z)$  the uniform branches of the above
multivalued function. By direct computation, we find as expected that $B\circ S_{j}=B,$ $j=0,1,...,n-1.$ If we conveniently choose the principal
branch of $(27)$ as in $(22),$ then we have again $\ S_{j}\circ
S_{k}=S_{k}\circ S_{j}=T_{j}\circ S_{k}=$ $S_{j}\circ T_{k}=S_{j+k}$ \ for
every $j$ and $k,$ and we can easily infer that the group generated by the
transformations $S_{j}$ and $z\rightarrow e^{2\pi i/n}z$ is the group of covering
transformations of $(\widehat{\mathbb C},B).$

\section{Infinite Blaschke Products}

\ 

It is known (see \cite{[2]} and \cite{[3]}) that if the Blaschke product $w=$ $B(z)$ has
the (generalized) Cantor set $E$ as the set of cluster points of its zeros,
then the equation $B(z)=1$ has infinitely many distinct solutions in every
"removed" open arc $I_{n} \subset \partial D\, \backslash \, E.$ The only
cluster points of these solutions are the end points of $I_{n}.$ Therefore, $\ds I_{n}=\cup _{k=-\infty }^{+\infty }\Gamma _{n,k},$ where $\Gamma _{n,k}$ are
half open arcs of the unit circle, which are mapped by $B$ continuously and
bijectively on the unit circle in the $w$-plane and they accumulate to the
ends of $I_{n}$. Every fundamental domain $\Omega _{n,k}$ of $B$ contains the
interior of a unique arc $\Gamma _{n,k}.$ These domains also accumulate  to
every point of $E.$ Moreover, if $K\subset \mathbb{C}
\, \backslash \, (A\cap E)$ is a compact set, then there is only a finite
number of domains $\Omega _{n,k}$ such that $\Omega _{n,k}\cap K\neq
\emptyset .$ Here $A=\{1/\overline{a}_{n,k},n\in N,k\in Z\}$ is the set
of the poles of $B.$\vspace{.1in}

We have seen in Section 2 that, for the type of Blaschke product studied
there, the fundamental domains corresponding to $a_{2}$ could be made as
small as we wanted if $a_{2}$ was moved close enough to the unit circle. We
conjecture here that this property is true for any Blaschke product.
Moreover, for an infinite Blaschke product $B$ with $E$ a Cantor set, if we
choose the sequence of partial products $\{B_{m}\}$ adding always $a_{k}$ with
non decreasing module, then there are sequences of fundamental domains $\{\Delta _{m,j}\}$ of $B_{m}$ convergent to every fundamental domain $%
\Omega _{n,k}$ of $B.$ In other words, not only $\{B_{m}\}$ approximates $B,$
but the fundamental domains of $\{B_{m}\}$ "approximate" every fundamental
domain of $B.$ Consequently, for $m$ large enough, the picture of fundamental
domains of $B_{m},$ as well as the mapping of those domains by $B_{m}$ give
a fairly accurate description of the first $m$ fundamental domains of $B,$
respectively of the mapping of those domains by $B.$ The remaining
fundamental domains of $B$ can be enclosed in the complement of a compact
set $K\subset \mathbb{C}
\, \backslash \, (A\cup E).$ More exactly, if $K$ is given, then there exists  $%
m_{0}\in 
\mathbb{N}
$ depending only on $K,$ such that for any $m\geq m_{0}$ only the first $%
m_{0}$ fundamental domains intersect $K.$ Then, to characterize $B,$ we can
simply deal with a $B_{m}$ for a large $m.$ Solving, as usual, the equations $%
B_{m}(z)$ $=1,$ and $B_{m}^{\prime }(z)=0$ is in principle a feasible task,
and so is the corresponding simultaneous continuation.
\vspace{.1in}

Let us consider the example of the infinite Blaschke product $ B$ defined by the
Blaschke sequence $\{a_{n}^{(k)}\},$ where
\begin{equation}a_{n}^{(k)}=\left[1-\frac{1}{(n+1)^{2}}\right]\omega _{k},\ \ \ k=0,1,2;\ \ n\in
\mathbb N,\label{eq30}\end{equation}

and $\omega _{k}$ are the roots of order three of  unity. If we let $%
a_{n}^{(0)}=a_{n},$ then

$$B(z)=\prod\limits_{n-1}^{\infty }\frac{z^{3}-a_{n}^{3}}{
a_{n}^{3}z^{3}-1}.$$
We notice that for this Blaschke product we have $E=\{\omega _{k},\  
k=0,1,2\}.$ We  illustrate by a picture the fact that $B$ is fairly
well approximated by the partial product having just nine factors, namely
the product obtained letting $n\in \{1,2,3\},$ $k\in \{0,1,2\}.$ An easy
computation shows that:
\begin{equation}B_{9}(z)=\frac{z^{9}-\sigma _{1}z^{6}+\sigma _{2}z^{3}-p}{pz^{9}-\sigma
_{2}z^{6}+\sigma _{1}z^{3}-1},\label{eq31}\end{equation}
where$$\sigma _{1}=a_{1}^{3}+a_{2}^{3}+a_{3}^{3},\qquad \sigma
_{2}=a_{1}^{3}a_{2}^{3}+a_{1}^{3}a_{3}^{3}+a_{2}^{3}a_{3}^{3},\qquad 
p=a_{1}^{3}a_{2}^{3}a_{3}^{3}.$$

\ \ 

The following equations are equivalent.
$$\begin{array}{l}B_{9}(z)=1\\ \ \\ (1-p)z^{9}+(\sigma _{2}-s_{1})z^{6}+(\sigma _{2}-s_{1})z^{3}+(1-p)=0\\ \ \\
(z^{3}+1)[(1-p)z^{6}-(1-p+\sigma _{1}-\sigma _{2})z^{3}+(1-p)]=0\end{array}$$
The last equation can be easily solved, and we obtain the solutions:
$$z_{k}^{(1)}=e^{i\pi /3}\omega _{k},\qquad z_{k}^{(2)}=e^{-i\theta
/3}\omega _{k},\qquad z_{k}^{(3)}=e^{i\pi /3}\omega _{k},$$
where $\omega _{k}$ are again the roots of order three of  unity and $%
\theta =11.5^{\circ}$. \ It is obvious that these nine solutions are all on
the unit circle, as we expected. Now we need the solutions of the equation $%
B_{9}^{\prime }(z)=0.$ We notice that $\ b_{0}=0$ is a double solution of
this equation. After the substitution $u=z^{3},$ the remaining solutions
are those of the equation $$ \frac{d}{du}\left(\frac{u^{3}-\sigma _{1}u^{2}+\sigma _{2}u-p}{pu^{3}-\sigma
_{2}u^{2}+\sigma _{1}u-1}\right)=0.$$
An elementary computation shows that this equation is equivalent to:
$$(s_{2}-p\sigma _{1})(u+1/u)^{2}+2(p\sigma _{2}-\sigma
_{1})(u+1/u)+(\sigma _{1}+p)^{2}+4(1-p^{2})-(1+\sigma _{2})^{2}=0$$ 

We notice that the images by $B_{9}$ of all these branch points are real,
hence we  expect to obtain fundamental domains for $B_{9}$ by
performing simultaneous continuation over the real axis of the $w$-plane
starting from the solutions of $B_{9}(z)=1.$ However, since the curves
obtained in this way pass through $a_{n}^{(k)},$ they cannot border
fundamental domains, given the fact that $a_{n}^{(k)}$ are simple zeros and
therefore they should remain inside of such domains. In order to avoid this
contradiction, we turn around the origin of the $w$-plane on a small
half-circle when performing simultaneous continuation.  Figure 6(a) below
shows the fundamental domains of $B_{9}$ obtained in this way. \vspace{.1in}

  \begin{center} \includegraphics[width=2.2truein,height=2.2truein]{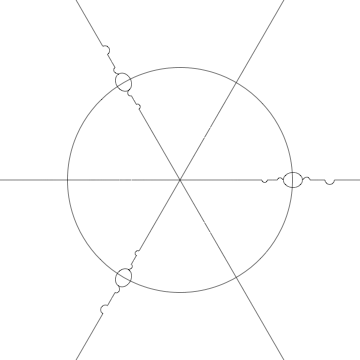} \\ Figure 6(a) \end{center} \vspace{.1in}
  
We perform  the same computation for $B_{12}$ and we obtain  Figure
6(b). 

 \vspace{.1in}

  \begin{center}  \includegraphics[width=2.2truein,height=2.2truein]{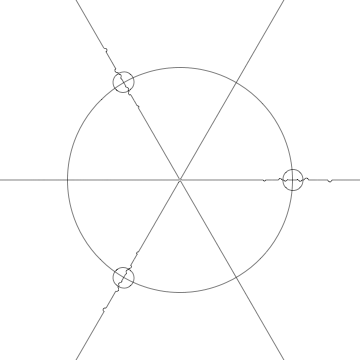} \\ Figure 6(b) \end{center} \vspace{.1in}

The unbounded fundamental domains remained practically the same, while
every bounded domain simply split into
two. If we focus just on a compact set $ K\subset 
\mathbb{C}
\, \backslash \, \cup _{k=0}^{3}D_{k},$ where $D_{k}$ are discs of radius $0.05$ centered
 at $\omega _{k}$ respectively,   then, starting with a given $%
m=9$, all these domains will belong to $\mathbb{C}
\, \backslash \, K$ and will be out of sight. In other words, as
long as we deal only with $K,$ the finite Blaschke product $B_{m}$ is the
same as the infinite Blaschke product $B.$\vspace{.1in}

Figure 5 illustrates the approximation of the infinite  Blaschke product defined by the sequence $(28)$ by the partial products $B_9$ and $B_{12}$.   Figures 5(g) and (h) zoom in  close to $z=1$ on the  fundamental domains of $B_9$ and $B_{12}$ respectively.  Figures 5(a) and (b) show the pre-images under $B_9$ and $B_{12}$ respectively of concentric annuli. Several of these annuli have been displayed in figures 5(e) and (f). The complete collection of annuli is available at \cite{cri}. Figures 5(c) and (d) zoom in close to $z=1$ on the respective pre-images given in figures 5(a) and (b). The pre-images under $B_9$ and $B_{12}$ 
 are hardly distinguishable, the
difference appearing only in the zoom, where in Figure 5(d) we can notice
a very small additional  fundamental domain.

\newpage
\begin{center} Figure 1
\end{center}

\vspace{.2in}

\noindent  \includegraphics[width=3.3truein,height=3.3truein]{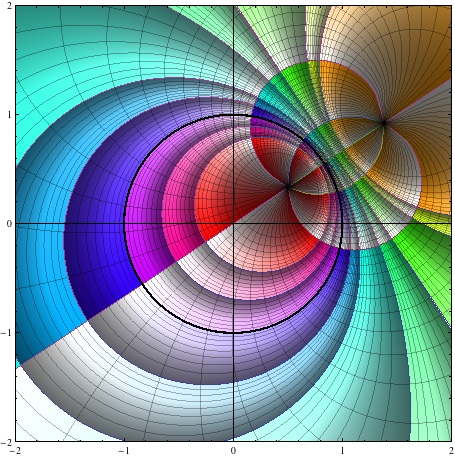} 
  \includegraphics[width=3.3truein,height=3.3truein]{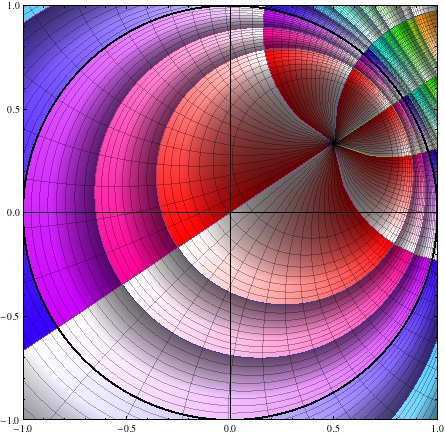}\newline
  \hspace*{1.5in} $(a)$ \hspace{3.1in} $(b)$\vspace{.2in}

\hspace*{2in}  \includegraphics[width=2.2truein,height=1.65truein]{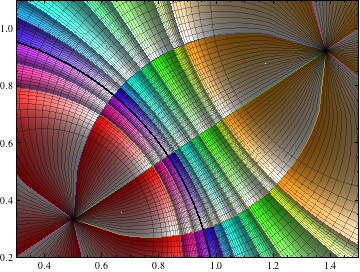} \newline

 \hspace*{3in}(c)

  \vspace{.2in}
  
   \noindent  \begin{tabular}[t]{lll}\includegraphics[width=2.1truein,height=2.1truein]{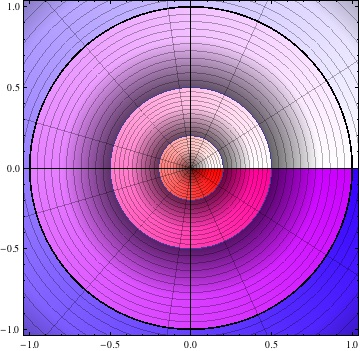} 
      \includegraphics[width=2.1truein,height=2.1truein]{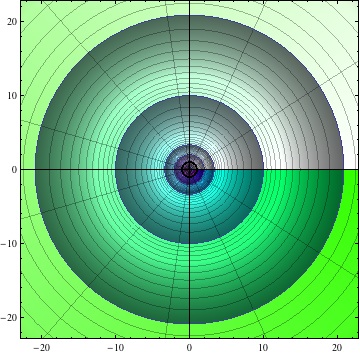} 
       \includegraphics[width=2.1truein,height=2.1truein]{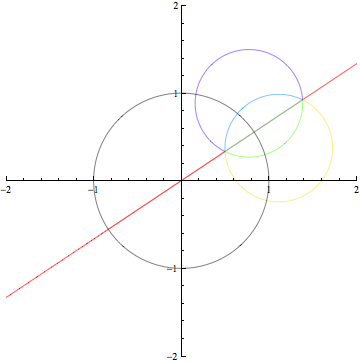}\end{tabular} 
 \newline
 \hspace*{.9in} $(d)$ \hspace{2in} $(e)$
   \hspace{1.7in} $(f)$\vspace{.1in}  

\newpage

\begin{center} Figure 2 
\end{center}\vspace{.2in}

\noindent  \includegraphics[width=3.3truein,height=3.3truein]{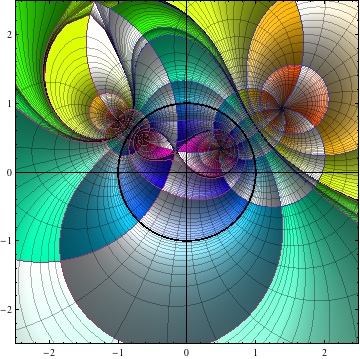} 
  \includegraphics[width=3.3truein,height=3.3truein]{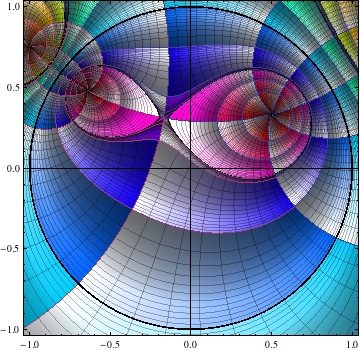} \newline
  \hspace*{1.5in} $(a)$ \hspace{3.1in} $(b)$\vspace{.2in} 

\noindent\hspace*{-.1in}  \begin{tabular}{lll}\includegraphics[width=2.3truein,height=1.7truein]{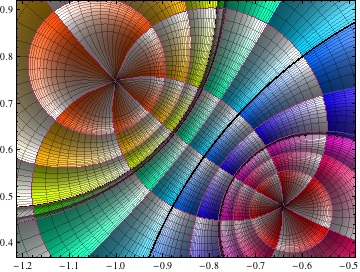} \includegraphics[width=2.1truein,height=2.1truein]{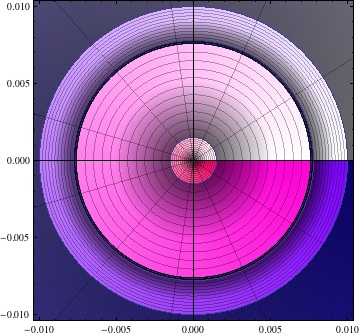} \includegraphics[width=2.1truein,height=2.1truein]{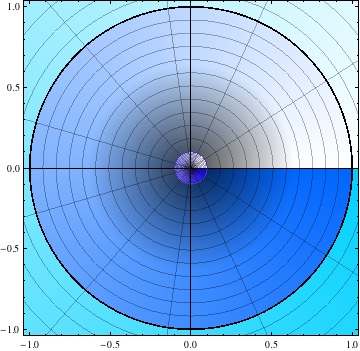} \end{tabular}

\hspace*{1in} (c) \hspace{2in} (d)\hspace{1.7in} (e)\vspace{.2in} 
  
   \noindent  \begin{tabular}[t]{lll}  \includegraphics[width=2.1truein,height=1.75truein]{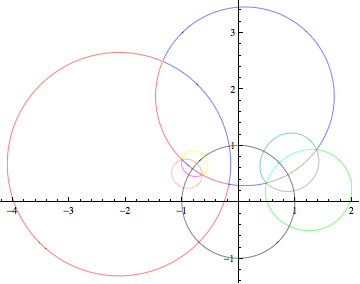}  \includegraphics[width=1.801truein,height=2.2truein]{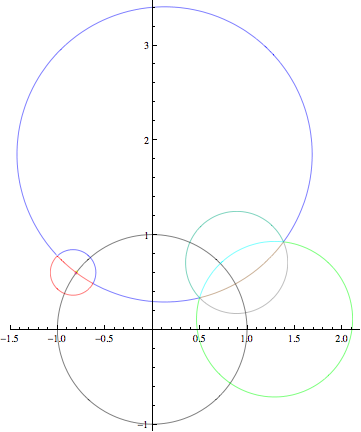} \includegraphics[width=1.801truein,height=2.2truein]{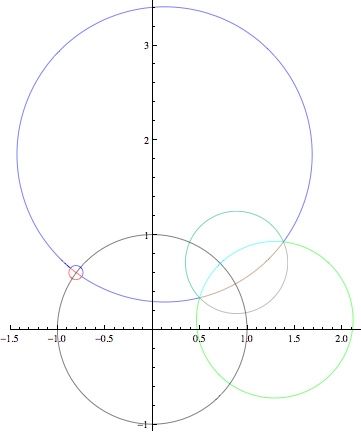}\end{tabular}
   \newline
  \hspace*{.9in} $(f)$ \hspace{1.7in} $(g)$
   \hspace{1.5in} $(h)$\vspace{.1in}

   \newpage
   
   \begin{center} Figure 3
\end{center}

\vspace{.2in}

\noindent  \includegraphics[width=3.3truein,height=3.3truein]{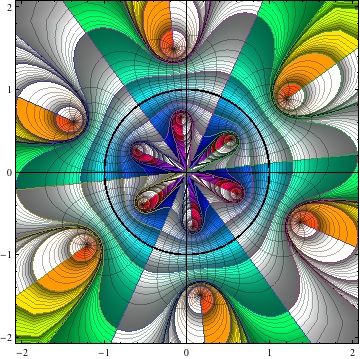} 
  \includegraphics[width=3.3truein,height=3.3truein]{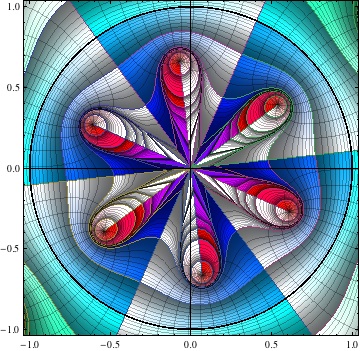} \newline
   \hspace*{1.5in} $(a)$ \hspace{3.1in} $(b)$\vspace{.2in}

\hspace*{.4in}\includegraphics[width=2.125truein,height=1.7truein]{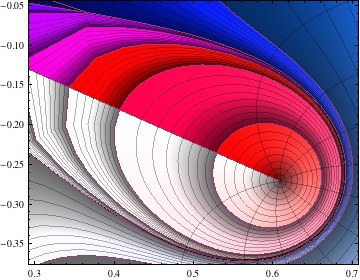} \hspace{1in}\includegraphics[width=2.2truein,height=2.2truein]{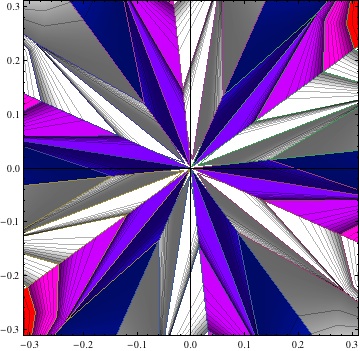}\newline
\hspace*{1.5in} (c)\hspace{3in}(d)\vspace{.2in}

   \noindent  \begin{tabular}[t]{lll}\includegraphics[width=2.1truein,height=2.1truein]{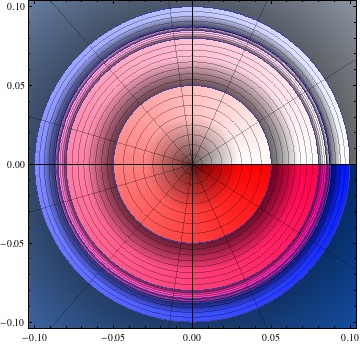} 
      \includegraphics[width=2.1truein,height=2.1truein]{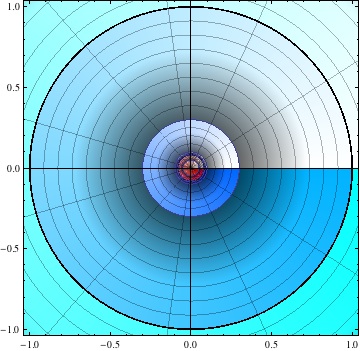} 
       \includegraphics[width=2.1truein,height=2.1truein]{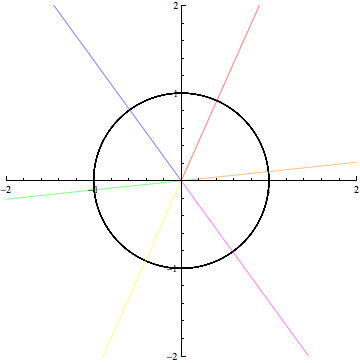} \end{tabular}
 \newline
    \newline
  \hspace*{.9in} $(e)$ \hspace{2in} $(f)$
   \hspace{1.7in} $(g)$\vspace{.1in} 
\newpage

\begin{center} Figure 4
\end{center}

\vspace{.2in}

\noindent  \includegraphics[width=3.3truein,height=3.3truein]{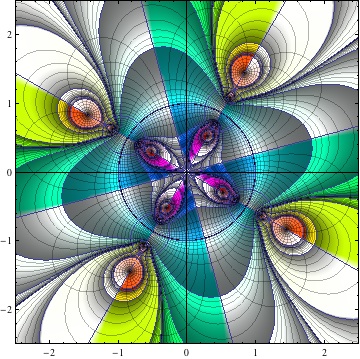} 
  \includegraphics[width=3.3truein,height=3.3truein]{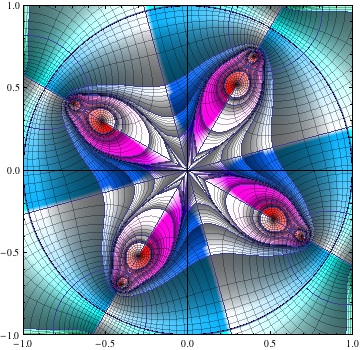}\newline
  \hspace*{1.5in} $(a)$ \hspace{3.1in} $(b)$\vspace{.5in}

  \noindent \hspace*{.5in} \includegraphics[width=2.2truein,height=1.6truein]{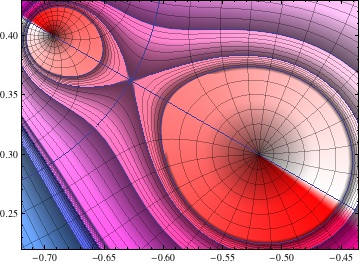} \hspace{.7in}
  \includegraphics[width=2.2truein,height=1.6truein]{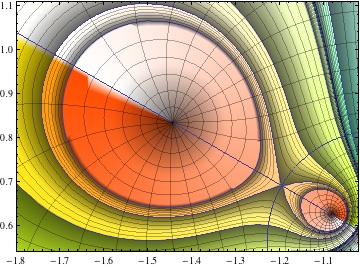}\newline
  \hspace*{1.5in} $(c)$ \hspace{2.7in} $(d)$\vspace{.2in}

   \noindent \begin{tabular}[t]{lll}\includegraphics[width=2.1truein,height=2.1truein]{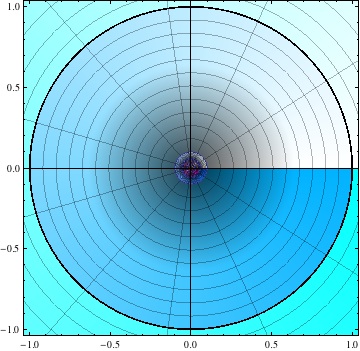} 
      \includegraphics[width=2.1truein,height=2.1truein]{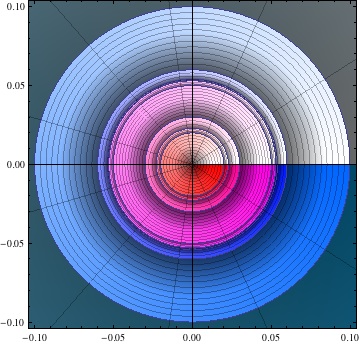} 
       \includegraphics[width=2.1truein,height=2.1truein]{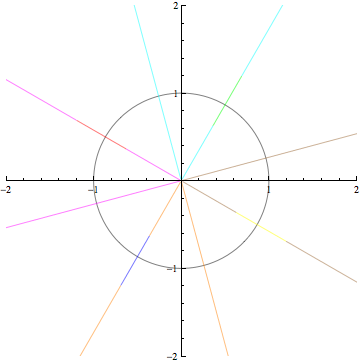}\end{tabular} 
 \newline
     \newline
  \hspace*{.9in} $(e)$ \hspace{2in} $(f)$
   \hspace{1.7in} $(g)$\vspace{.1in}    
   \newpage
   
\begin{center} Figure $5$\end{center} \vspace{.2in}

\noindent  \includegraphics[width=3.3truein,height=3.3truein]{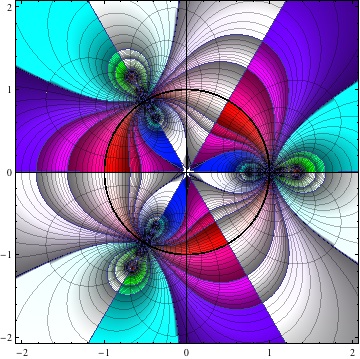} 
  \includegraphics[width=3.3truein,height=3.3truein]{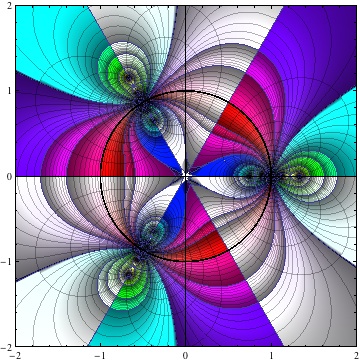}\newline
  \hspace*{1.5in} $(a)$ \hspace{3.1in} $(b)$\vspace{.5in}

  \noindent \hspace*{.5in} \includegraphics[width=2.4truein,height=1.2truein]{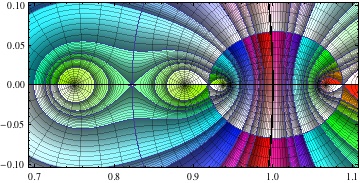} \hspace{.7in}
  \includegraphics[width=2.4truein,height=1.2truein]{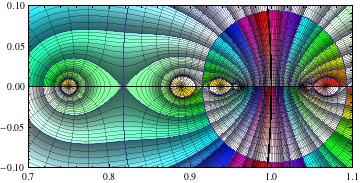}\newline
  \hspace*{1.7in} $(c)$ \hspace{3in} $(d)$\vspace{.2in}

 \noindent \hspace*{-.2in} \begin{tabular}[t]{lll}{  \includegraphics[width=2.1truein,height=2.1truein]{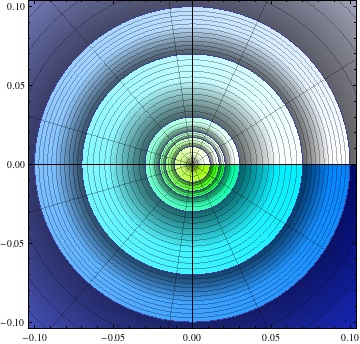} 
      \includegraphics[width=2.1truein,height=2.1truein]{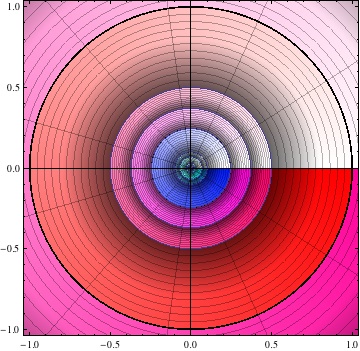}
      \begin{tabular}[b]{c} \includegraphics[width=2.125truein,height=.5truein]{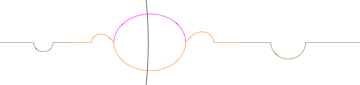}\\ $(g)$\\\mbox{\ }\\ \mbox{\ }\\
      \includegraphics[width=2.125truein,height=.5truein]{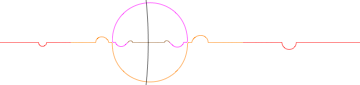} \\ $(h)$\\ \mbox{\ }\end{tabular}}\ \end{tabular}

  \hspace*{.8in}  $(e)$ \hspace{1.9in} $(f)$ \vspace{.2in}
  
\newpage

\section{Technical details} 

All images for this article have been created using the software  \emph{Mathematica 6} on a MacBookPro with a 2.33 GHz Intel Core 2 Duo processor. Sample code is available on the website of the project \cite{cri}. 

\section*{Acknowledgements}

The authors would like to thank Maxim Rytin for a very helpful hint on Mathematica and Cristi Rinklin for artistic advice.

\end{document}